\documentclass{amsart} 
\usepackage{graphics}
\usepackage{color}
\usepackage{epic}

\newtheorem{theo}{\bf Theorem}[section]
\newtheorem{propo}[theo]{\bf Proposition}
\newtheorem{lemma}[theo]{\bf Lemma}
\newtheorem{conj}[theo]{\bf Conjecture}

\newtheorem{coro}[theo]{\bf Corollary}

\newcommand{\cE}{{\mathcal E}}
\newcommand{\cV}{{\mathcal V}}
\newcommand{\ka}{{\kappa}}
\newcommand{\rE}{{\rm E}}
\newcommand{\Prob}{{\rm P}}

\begin{document}

\title[Exact expectations for random graphs and assignments]
{Exact expectations for random graphs and assignments}
\author{Henrik Eriksson}
\address{NADA \\
   KTH \\
   SE-100 44 Stockholm, Sweden}
\email{henrik@nada.kth.se}
\author{Kimmo Eriksson}
\address{IMA, M{\"a}lardalens h{\"o}gskola \\ 
   Box 883 \\
   SE-721 23 V{\"a}ster{\aa}s, Sweden}
\email{kimmo.eriksson@mdh.se}
\author{Jonas Sj{\"o}strand}
\address{NADA \\ 
   KTH \\
   SE-100 44 Stockholm, Sweden}
\email{jonass@nada.kth.se}
\date{November 12, 2000}

\begin{abstract}
For a random graph on $n$ vertices where the edges appear with individual
rates, we give exact formulas for the expected time at which the number 
of components has gone down to $k$ and the expected length of the 
corresponding minimal spanning forest.

For a random bipartite graph we give a formula for the expected time at 
which a $k$-assignment appears. This result has bearing upon the random 
assignment problem.

\end{abstract}
\maketitle

\section{Introduction}\noindent
The study of random graphs, initiated by Erd\H{o}s forty years ago, has
been preoccupied by asymptotic results. In this paper, we will derive
exact enumerative formulas in some cases where previously only asymptotics 
were known. 

By a random graph $G(t)$ on $n$ vertices we mean a stochastic 
process starting with an edgeless graph at $t=0$ and
where edge after edge appears at random times. We take these times
to be independent exponentially distributed variables,
so the probability of edge $e_{ij}$ not having appeared at time $t$
is $\exp(-a_{ij}t)$. The rates $a_{ij}$ may be viewed as labelling
the edges of $K_n$ and any nonnegative edge labelling specifies 
such a random graph process.

The number of graph components at time $t$ is a thoroughly studied
statistic and asymptotics for its expectation and distribution are
known in many cases (see Janson \cite{SJ}). Of particular 
importance is the time $T_n$ at which the graph becomes one component,
i.e.\@ gets a spanning tree, and also the minimal length $L_n$ of such a
tree, i.e.\@ the sum of its edge times. 
(By Kruskal's algorithm the first tree has indeed minimal length.)
For the simplest case when all
rates are 1, there is the following beautiful result by Frieze.
$$\lim_{n\rightarrow
\infty} \rE(L_n) = \zeta(3)=1+\frac{1}{2^3}+\frac{1}{3^3}+\frac{1}{4^3}+\dots
= 1.202\dots$$
The first published version of the above
result seems to be by Fenner and Frieze \cite{FF}, where
the distribution used was the uniform distribution on $[0,1]$.
Later Frieze
\cite{F} showed that the same result holds
for any distribution function $F$ such that
$F'(0)$ exists
and equals $1$. The Exp(1) distribution is an obvious example.

In the first theorem of this paper, we will give simple
exact formulas for the following more general statistics.
\begin{eqnarray*}
\rE(T_n(k))&=&\rE(t\, |\, G(t)\, \mbox{becomes a $k$-component graph}) \\
\rE(L_n(k))&=&\rE(\mbox{length of the minimal spanning $k$-forest})
\end{eqnarray*}

\noindent
Note that our formulas work for arbitrary rates  $a_{ij}$.
The special case where all  $a_{ij}=1$ turns out to be a sum over
all integer partitions of $n$ (Corollary \ref{co:min_tree}).
In the general case we will use the following notation.
Let the rates  $a_{ij}$ label the edges $\cE_n$ of $K_n$.
For an edge $e_{ij}\in\cE_n$, let $[e_{ij}]=a_{ij}$ and for a subset 
$A\subseteq \cE_n$
of edges, let $[A]$ denote the rate sum $\sum_{e_{ij}\in
A} a_{ij}$. A subset $A\subseteq \cE_n$ is a {\em clique} if it spans a
complete subgraph. Any partition of the vertices of $K_n$ into $j$ parts
defines a subgraph consisting of $j$ cliques. If $B$ is the edge set
of such a subgraph, we write $B\in \mbox{Clique}_j$. 

\begin{theo}\label{th:min_forest}
For a random graph $G(t)$ on $n$ vertices with rates $a_{ij}$
\begin{eqnarray*}
\rE(T_n(k))&=&  \sum_{j=k+1}^n \sum_{B\in \mbox{\rm\footnotesize Clique}_j} 
\frac{\tau (j,k)}{[\cE_n-B]}\\
\rE(L_n(k))&=& \sum_{j=k+1}^n 
\sum_{B\in \mbox{\rm\footnotesize Clique}_j} 
\frac{\lambda (j,k)}{[\cE_n-B]}
\end{eqnarray*}
where  $\tau (j,k)=\sum_{i=k+1}^j s(j,i)$ and
$\lambda (j,k)=\sum_{i=k+1}^j (i-k)s(j,i)$ are sums of
signed Stirling numbers of the first kind.
\end{theo}

\begin{table}[htb] \label{theTable}
\begin{center}
\begin{tabular}{|r|rrrrrr|}
\hline 
$s$ & 0 & 1 & 2 & 3 & 4 & 5\cr 
\hline 
0   &  1 &    &    &    &    &  \cr
1   &  0 &  1 &    &    &    &  \cr
2   &  0 & -1 &  1 &    &    &  \cr
3   &  0 &  2 & -3 &  1 &    &  \cr
4   &  0 & -6 & 11 & -6 &  1 &  \cr
5   &  0 & 24 &-50 & 35 &-10 & 1\cr
\hline
\end{tabular}
\ \ \ 
\begin{tabular}{|r|rrrrr|}
\hline 
$\tau$ & 0 & 1 & 2 & 3 & 4 \cr 
\hline 
0   &  0 &    &    &    &  \cr
1   &  1 &    &    &    &  \cr
2   &  0 &  1 &    &    &  \cr
3   &  0 & -2 &  1 &    &  \cr
4   &  0 &  6 & -5 &  1 &  \cr
5   &  0 &-24 & 26 & -9 & 1\cr
\hline
\end{tabular}
\ \ \ 
\begin{tabular}{|r|rrrrr|}
\hline 
$\lambda$ & 0 & 1 & 2 & 3 & 4 \cr 
\hline 
0   &  0 &    &    &    &  \cr
1   &  1 &    &    &    &  \cr
2   &  1 &  1 &    &    &  \cr
3   & -1 & -1 &  1 &    &  \cr
4   &  2 &  2 & -4 &  1 &  \cr
5   & -6 & -6 & 18 & -8 & 1\cr
\hline
\end{tabular}

\end{center}
\caption{The recursion $s(j\!+\!1,i)=s(j,i\!-\!1)-j\, s(j,i)$ is shared by 
$\tau$ and $\lambda$}
\end{table}

\noindent
For example, $K_3$ has three subgraphs in $\mbox{Clique}_2$ 
and one in $\mbox{Clique}_3$, so  
$$ \rE(T_3(1))=\frac{1}{a_{12}+ a_{13}}+\frac{1}{a_{12}+ a_{23}}+
\frac{1}{a_{13}+ a_{23}}- \frac{2}{a_{12}+ a_{13}+a_{23}}.
$$

\noindent
If we put all $a_{ij}=1$ and also $k=1$, we should get an expression
that tends to Frieze's result when $n$ goes to
infinity.
We use the
notation 
$\binom{n}{{1^{e_1},2^{e_2},\dots}}$ for the multinomial number
$\frac{n!}{{(1!)^{e_1}(2!)^{e_2}\dots}}$.

\begin{coro}\label{co:min_tree}
For a random graph $G(t)$ on $n$ vertices with rates $1$
\begin{eqnarray*}
\rE(T_n(1))=&\sum {\binom{n}{ {1^{e_1},2^{e_2},\dots}}}
{\binom{e}{{e_1,e_2,\dots}}}
        \frac{(-1)^e}{{\frac{e}{2}}
(n^2-(1^2e_1 + 2^2e_2 + \dots))}\\
\rE(L_n(1))=&\sum {\binom{n}{{1^{e_1},2^{e_2},\dots}}}
{\binom{e}{{e_1,e_2,\dots}}}
        \frac{(-1)^e}{{\binom{e}{2}}
(n^2-(1^2e_1 + 2^2e_2 + \dots))}
\end{eqnarray*}
where the sums are over integer partitions of $n$ into two or more parts, 
i.e.\@ nonnegative integer sequences 
$(e_i)$ such that 
$e=e_1+e_2+\dots \ge 2$ and 
$1e_1 + 2e_2 + \dots = n$.
\end{coro}

\noindent
For example, $ \rE(L_3(1))= {\binom{3}{ 1,2}}{\binom{2}{1,1}}/
{\binom{2}{2}}(9\!-\!1\!-\!4)
-{\binom{3}{1,1,1}}{\binom{3}{3}}/{\binom{3}{2}}(9\!-\!3) = 
{\frac{7}{6}}$.
Clearly, the terms of the formulas grow exponentially with
$n$. It is an intriguing consequence of Frieze's theorem
that, nevertheless, the last sum must tend to $\zeta(3)$ as
$n\rightarrow \infty$.

\subsection{Random assignment}
Our reason for studying this
problem was to get a grip on the {\em random assignment
problem} which is the study of minimal assignments in a
random bipartite graph $G(t)\subseteq K_{m,n}$.
As before, the time of appearance of edge $e_{ij}$ is
assumed to be exponentially distributed with
rate $a_{ij}$ and we consider these rates as labelling
the edges of the complete bipartite graph $K_{m,n}$.

A $k$-assignment is a disjoint set of $k$ edges. The following 
two statistics are of special interest.
\begin{eqnarray*}
\rE(T_{mn}(k))&=&\rE(t\, |\, G(t)\, \mbox{obtains a $k$-assignment}) \\
\rE(L_{mn}(k))&=&\rE(\mbox{length of the minimal $k$-assignment})
\end{eqnarray*}
The analogy with the previous definitions is deceptive, for here the
$k$-assignment of $T_{mn}(k)$ is not the $k$-assignment of $L_{mn}(k)$
in general! Since the length of an assignment is the sum of the 
appearance times of the edges involved, the first $k$-assignment 
to be completed need not have the minimal length. Because of this,
the computation of $\rE(L_{mn}(k))$ seems to be quite difficult
and we have obtained formulas only for the simplest cases.

For $k\!=\!m\!=\!n$, we have the original random assignment problem.
The following conjecture, proposed by Parisi \cite{P} in 1998,
has been verified up to $n=7$.
\begin{conj}
  $\rE(L_{nn}(n))=1+1/4+1/9+\dots+1/n^2$.
\end{conj}

\noindent
Note that this expression tends to $\zeta(2)$ as $n\rightarrow 
\infty$. This asymptotic conjecture, stated by M\'ezard and 
Parisi in 1985 \cite{MP}, was proved by Aldous \cite{Aldous}.

Generalizations of the conjecture to general $k,m,n$ have been
proposed by Linusson and W\"astlund \cite{LW} and by Buck, Chan
and Robbins \cite{BCR}, who also consider more general rates
$a_{ij}$. As far as we know, however, formulas for arbitrary rates have not
been obtained before.  

For
the expected time until the first $k$-assignment appears, our
method gives the following formula.

\begin{theo}\label{th:time_assignment}
$$\rE(T_{mn}(k))=\sum_{B\in \mbox{\rm\footnotesize Tabloid}_k} 
\frac{S_{m,n,k}(B)}{[\cE_{m,n}-B]}$$
\end{theo}

\noindent The numerators $S_{m,n,k}(B)$ are certain products of binomial
coefficients to be defined in Section 3. A proper edge subset $B$ is called tabloidal
if the matrix entries $e_{ij}\in B$ fill a tableau shape (Young diagram),
possibly after permutation of rows and columns. It is in 
$\mbox{\rm Tabloid}_k$ if the inner corners of this
tableau satisfy $i+j<k$ and the outer corners $i+j\ge k$.

\definecolor{light}{gray}{.55}
\begin{figure}
\begin{picture}(60,60)(0,0)
\put(10,10){\framebox(40,30)}
\dashline[-7]{2}(10,20)(30,40)
\end{picture}
\begin{picture}(60,60)(0,0)
\put(10,10){\framebox(40,30)}
\put(10,33){\colorbox{light}{\makebox(4,4){}}}
\dashline[-7]{2}(10,20)(30,40)
\end{picture}
\begin{picture}(60,60)(0,0)
\put(10,10){\framebox(40,30)}
\put(10,33){\colorbox{light}{\makebox(34,4){}}}
\dashline[-7]{2}(10,20)(30,40)
\end{picture}
\begin{picture}(60,60)(0,0)
\put(10,10){\framebox(40,30)}
\put(10,13){\colorbox{light}{\makebox(4,24){}}}
\dashline[-7]{2}(10,20)(30,40)
\end{picture}
\caption{The four tableau shapes in $\mbox{\rm Tabloid}_2$.}
\end{figure}
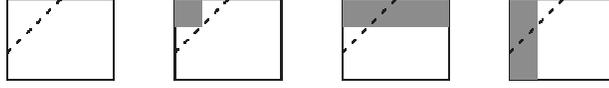

With $k=2$, for example, four tableau shapes are involved
and the expression becomes
$$\rE(T_{mn}(2))=
\frac{\binom{m-1}{1}\binom{n-1}{1}}{[\cE_{m,n}]}
-\sum_{i,j}\frac{1}{[\cE_{m,n}\!-\!e_{ij}]}
+\sum_{i}\frac{1}{[\cE_{m,n}\!-\!\mbox{row}_i]}
+\sum_{j}\frac{1}{[\cE_{m,n}\!-\!\mbox{col}_j]}.$$

\noindent
If we try to compute $L_{mn}(k)$ as
$T_{mn}(1)+T_{mn}(2)+\cdots+T_{mn}(k)$
we get a value that is too small, since the subassignments of the assignment of
minimal length need not be minimal for $T_{mn}$.
For small $k$ we have been able to calculate $\rE(L_{mn}(k))$, for example
$$\rE(L_{mn}(2))=\rE(T_{mn}(1))+\rE(T_{mn}(2))+\sum_{i,j}\frac{a_{ij}[\mbox{row}_i-e_{ij}]
[\mbox{col}_j-e_{ij}]}{[\cE_{m,n}][\cE_{m,n}-e_{ij}][\cE_{m,n}-\mbox{row}_i][\cE_{m,n}-\mbox{col}_j]}.$$
The derivation of this expression and similar ones will appear in our forthcoming
paper \cite{EES}.

\section{Minimal spanning $k$-forests}\label{sec:forest}
\noindent
This section is devoted to the proof of Theorem \ref{th:min_forest}. The
formula to be proved is of the form
$\rE(T_n(k))=  \sum c_B/[\cE_n-B]$
summed over all edge sets $B$ (the constant $c_B$ being zero for
most $B$). Now, the terms on the right-hand side
can be interpreted according to an elementary property of exponentially
distributed variables: the expected minimum of a set of such variables
is the reciprocal of the total rate. In our case, this means
$\rE(\mbox{minimal appearance time of an edge in } A)=1/{[A]}$.

In view of this, we are going to prove a derandomized version
of the formula, obtained by removing the expectation operator.
\begin{equation}
\label{eq:derandomized}
T_n(k)=  \sum_B c_B\ \min(\cE_n-B)
\end{equation}
\noindent
In this version, $T_n(k)$ is no longer a stochastic variable. The graph
$G(t)$ may have been constructed by any process, random or nonrandom,
with edge $e_{ij}$ appearing at time $t_{ij}$ and by $\min(\cE_n-B)$ we
mean $\min(t_{ij}\, |\, e_{ij}\notin B)$. Thus we have the 
problem of proving a purely combinatorial formula valid for any 
set of numbers $t_{ij}$.

In order to compute the time $T_n(k)$, we will 
study the Boolean lattice
of all subsets of edges. Each subset $A\subseteq \cE_n$
of edges determines a subgraph $(\cV_n,A)$ of $K_n$, where all $n$ original
vertices are still present. Let $\ka(A)$
 denote the number of connected components of this subgraph.

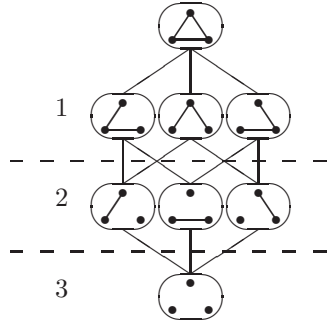
\begin{figure}[h]
\begin{picture}(60,110)(-40,-5)
\setlength{\unitlength}{0.6mm}
\put(0,0){\oval(14,10)}
\put(0,3){\circle*{2}}
\put(-4,-3){\circle*{2}}
\put(4,-3){\circle*{2}}

\put(0,20){\oval(14,10)}
\put(0,23){\circle*{2}}
\put(-4,17){\circle*{2}}
\put(4,17){\circle*{2}}
\put(-4,17){\line(1,0){8}}

\put(-15,20){\oval(14,10)}
\put(-15,23){\circle*{2}}
\put(-19,17){\circle*{2}}
\put(-11,17){\circle*{2}}
\put(-15,23){\line(-2,-3){4}}

\put(15,20){\oval(14,10)}
\put(15,23){\circle*{2}}
\put(19,17){\circle*{2}}
\put(11,17){\circle*{2}}
\put(15,23){\line(2,-3){4}}

\put(0,40){\oval(14,10)}
\put(0,43){\circle*{2}}
\put(-4,37){\circle*{2}}
\put(4,37){\circle*{2}}
\put(0,43){\line(2,-3){4}}
\put(0,43){\line(-2,-3){4}}

\put(-15,40){\oval(14,10)}
\put(-15,43){\circle*{2}}
\put(-19,37){\circle*{2}}
\put(-11,37){\circle*{2}}
\put(-15,43){\line(-2,-3){4}}
\put(-19,37){\line(1,0){8}}

\put(15,40){\oval(14,10)}
\put(15,43){\circle*{2}}
\put(19,37){\circle*{2}}
\put(11,37){\circle*{2}}
\put(15,43){\line(2,-3){4}}
\put(19,37){\line(-1,0){8}}

\put(0,60){\oval(14,10)}
\put(0,63){\circle*{2}}
\put(-4,57){\circle*{2}}
\put(4,57){\circle*{2}}
\put(0,63){\line(2,-3){4}}
\put(0,63){\line(-2,-3){4}}
\put(-4,57){\line(1,0){8}}

\put(0,5){\line(3,2){15}}
\put(0,5){\line(0,1){10}}
\put(0,5){\line(-3,2){15}}

\put(0,25){\line(3,2){15}}
\put(0,25){\line(-3,2){15}}
\put(0,35){\line(3,-2){15}}
\put(0,35){\line(-3,-2){15}}
\put(-15,35){\line(0,-1){10}}
\put(15,35){\line(0,-1){10}}

\put(0,55){\line(-3,-2){15}}
\put(0,55){\line(3,-2){15}}
\put(0,55){\line(0,-1){10}}

\dashline{3}(-40,10)(30,10)
\dashline{3}(-40,30)(30,30)

\put(-30,0){3}
\put(-30,20){2}
\put(-30,40){1}

\end{picture}
\caption{The Boolean subgraph lattice of $K_3$ partitioned according to the number
of components.}
\end{figure}

When
edges appear, one at a time, we follow a path upwards in
the Boolean lattice from the empty set to the complete edge
set. At some time $T_n(k)$, the path reaches a set $A$ where $\ka(A)=k$.
The first part of the following lemma simply states that this time 
is the sum of the time intervals spent in each node on the upward path.

\begin{lemma}\label{lm:time}
\begin{equation}\label{eq:lemsum}
T_n(k) = \sum_{\ka(A_i)>k}  T(A_i),
\end{equation}
where $A_i$ is the set consisting of the first $i$ edges (so 
$A_0=\emptyset$) and $T(A_i)$ is the time spent in configuration $A_i$ 
until another edge appears.
For any edge set $A$, the expression
\begin{equation}\label{eq:PE}
 T(A) =
  \sum_{B\subseteq
A} (-1)^{|A-B|}\ \min(\cE_n-B)
\end{equation}
evaluates to $T(A_i)$ if $A$ is an $A_i$ and to zero otherwise. 
\end{lemma}
\begin{proof}

Let $i=|A|$,
let $e_i$ be the last edge in $A$ to appear and let $t_i$ be
its time of appearance. For every $B\subseteq A$, let $\bar B$ be the
result of toggling $e_i$, that is either including or deleting it.
In most cases, the contributions of $B$ and $\bar B$ in (\ref{eq:PE})
will cancel, for $t_i$ will seldom be a minimum. The only case
for which $ \min(\cE_n-B)=t_i$ is when all earlier edges are in $B$.
But then $\bar B = A$ consists of the first $i$ edges, that is $A=A_i$,
and we have $\min(\cE_n-\bar B)-\min(\cE_n-B)=t_{i+1}-t_i=T(A_i)$.
Otherwise, all terms are cancelled by the toggling.
\end{proof}

\noindent
Combining the two formulas in Lemma \ref{lm:time} we obtain
the following
expression for the time at which the number of components of the graph $G(t)$ 
has gone down to $k$:
\begin{equation}\label{eq:sumB}
T_n(k) =  \sum_{\ka(A)>k} \sum_{B\subseteq A} 
(-1)^{|A-B|}\ \min(\cE_n-B)
=
  \sum_{B\subseteq \cE_n} \min(\cE_n-B) \sum_{B\subseteq
A,\;\ka(A)>k} (-1)^{|A-B|} 
\end{equation}

\noindent
We next observe that the last sum is zero for
many sets $B$. For suppose that $B$ has some component
that is not a clique (a complete subgraph). In other words,
there exists some edge $e$ that we can add to $B$ without
changing the number of components of $B$. This property
is obviously shared also by all supersets $A$ of $B$ where
$e$ is not already present. Pair off every such $A$ with $A+e$.
The pair's contribution to the sum is  $(-1)^{|A-B|}+(-1)^{|A+e-B|}=0$.
The entire sum is over a disjoint collection of such pairs,
so it is zero too.

We have now reduced the sum to sets
$B$ where all components are cliques.
If there are $j$ components, we write $B\in \mbox{Clique}_j$. 
It turns out that each clique can be contracted into a node of $K_j=(\cV_j,\cE_j)$:

\begin{lemma}\label{lm:j_comp}
If
$B\in \mbox{\rm Clique}_j$, the following equality of generating functions holds:
\begin{equation*}\label{eq:sumA}
 \sum_{B\subseteq A \subseteq \cE_n} (-1)^{|A-B|}x^{\ka(A)}
=
  \sum_{A'\subseteq \cE_{j}} (-1)^{|A'|}x^{\ka(A')}.
\end{equation*}
\end{lemma}
\begin{proof}
Fix $j$ and argue by induction over $n$. The equality is
trivially true for $n=j$ (in which case $B$ must be the
empty set). Suppose it holds for $n=i\ge j$. Now study $n=i+1$.
Some clique has at least two vertices, say $x$ and $y$.
We shall see that omitting vertex $y$ does not affect the sum.

If $A$ is an extension of $B$ including an edge from $y$ 
to a vertex $z$ outside the clique, then toggle edge 
$xz$ to obtain $\bar{A}$ (that is, the only difference between 
$A$ and $\bar{A}$ is whether the edge
$xz$ is included or not). Obviously $\ka(A)=\ka(\bar{A})$ so
they both contribute to the same term in the generating
function. Their joint contribution is zero. Hence, the existence of vertex
$y$ can be disregarded in the sum over edge sets $A$ that are extensions of 
$B$, so we are down to $i$ vertices again.
\end{proof}

The right-hand generating function
of (\ref{eq:sumA}) can be computed explicitly. It can be
derived as a special case of the prime example of the exponential
method in Stanley's book
\cite{St2}, but we prefer to give
a direct proof.

\begin{lemma}\label{lm:Stirling}
$$\sum_{A'\subseteq
\cE_{j}} (-1)^{|A'|}x^{\ka(A')}=x(x-1)(x-2)\cdots(x-j+1).$$
(This is the generating function $\sum_{i=1}^j s(j,i)x^i$
for the signed Stirling numbers of the first kind.)
\end{lemma}
\begin{proof} We argue by induction on $j$. The formula is true for $j=1$.
For $j>1$ vertices, all $A'$ such that at least two edges
involve
the $j$th vertex can be combined into zero-contributing
pairs by toggling the edge forming a triangle with the lexically
first pair of such edges. Hence, we need only consider $A'$
where the $j$th vertex is involved in zero edges (one possibility,
decreasing the number of components by one if removed),
or one edge ($j-1$ possibilities, changing sign if removed),
thus resulting in the factor $(x-j+1)$ times the generating
function for $j-1$ vertices.
\end{proof}

By the previous two
lemmas we can compute the sum for $B\in \mbox{Clique}_j$ as a sum
of Stirling numbers:
$$  \sum_{B\subseteq A \subseteq \cE_n,\;
\ka(A)>k} (-1)^{|A-B|} = 
    \sum_{A'\subseteq \cE_{j},\;
\ka(A')>k} (-1)^{|A'|} = 
    \sum_{i=k+1}^j s(j,i).
$$
We
can now plug this into (\ref{eq:sumB}) and obtain 
$$
T_n(k) = \sum_{j=k+1}^n
\sum_{i=k+1}^j s(j,i) \sum_{B\in \mbox{\footnotesize Clique}_j} 
 \min(\cE_n-B),
$$
thereby proving the first formula of Theorem \ref{th:min_forest}.

\subsection{Special cases}
The spanning tree special case (i.e.\@ $k=1$) can be simplified further.
Since the sum of all Stirling numbers $s(j,i)$ for a fixed
$j>1$ is zero, we have
$$
\sum_{i=2}^j s(j,i) = -s(j,1) = (-1)^j (j-1)!.
$$
(Recall that the signed Stirling numbers
of the first kind have the combinatorial
interpretation
that $s(j,i)$ is $(-1)^{j+i}$ times the number of permutations in $S_j$
with $i$ cycles.)

Specializing to all rates
$a_{ij}$ equaling 1, every partition $B$ of $K_n$ into
cliques of given sizes contributes equally. With $e_1$ cliques
of size 1,
$e_2$ cliques of size 2, etc., the number of
possibilities for $B$ is
$$
  \frac{n!}{{(1!)^{e_1}e_1!(2!)^{e_2}e_2!\dots}},
$$
and
the rate reciprocal $1/[\cE_n-B]$ is $2/(n^2-(1^2e_1
+ 2^2e_2 + \dots))$.  Let $e=e_1+e_2+\dots$
denote the number
of cliques. We must multiply everything by $(-1)^e(e-1)!=(-1)^e
e!/e$. Thus we obtain
$$\sum {\binom{n}{ {1^{e_1},2^{e_2},\dots}}}{\binom{e}{{e_1,e_2,\dots}}}
        \frac{(-1)^e}{\frac{e}{2}(n^2-(1^2e_1 + 2^2e_2 + \dots))},
$$
where the sum is
taken over all $e_1,e_2,\dots$ such that $1e_1+2e_2+\dots=n$.
This proves the first statement of Corollary \ref{co:min_tree}.

\subsection{Expected length of the minimal spanning $k$-forest}
By Kruskal's algorithm for spanning trees,
the minimal spanning $k$-forest in $K_n$ consists of the longest
edges of the first spanning $r$-forests for $r=n,n-1,n-2,\dots,k$.

Thus, the expected
length of the minimal spanning $k$-forest in $K_n$ is 
$$
\rE(L_n(k)) = \sum_{r=k}^n \rE(T_n(r))=
\sum_{r=k}^n\sum_{j=r+1}^n\sum_{i=r+1}^j s(j,i) \sum_{B\in \mbox{\footnotesize Clique}_j} 
\frac{1}{[\cE_n-B]}
$$
$$
=\sum_{j=k+1}^n\sum_{i=k+1}^j\sum_{r=k}^{i-1} s(j,i) \sum_{B\in \mbox{\footnotesize Clique}_j} 
 \frac{1}{[\cE_n-B]}
$$
proving the second half of Theorem \ref{th:min_forest}.
The second half of Corollary \ref{co:min_tree} is proved
in the same way as the first half.

\section{Expected time to a $k$-assignment in $K_{m,n}$}
\noindent
In this section we study a random bipartite graph 
$G(t)\subseteq K_{m,n}$. Instead of 
trees, we are waiting for a $k$-assignment (a collection of $k$ disjoint edges)
to appear. Let $\mu(A)$ denote the size of the 
maximal subset of disjoint edges in an edge set $A\subseteq \cE_{m,n}$. 

As edges appear in $G(t)$ we follow a path from the empty set and
upwards in the Boolean lattice, until we reach a set $A$ where $\mu(A)=k$.
Let $T_{mn}(k)$ be the waiting time. Precisely as in the previous section,
we have

\begin{equation}\label{eq:ass-lemsum}
 \rE(T_{mn}(k)) = \sum_{\mu(A)<k}\Prob(A)\rE(T_A).
\end{equation}
\noindent
As before, we have $\rE(T_A) = 1/[\cE_{m,n}-A]$ and obtain
\begin{equation}\label{eq:ass-PE}
 \Prob(A)\rE(T_A)=
  \sum_{B\subseteq A} \frac{(-1)^{|A-B|}}{[\cE_{m,n}-B]},
\end{equation}
and finally
\begin{equation}\label{eq:ass-sumB}
 \rE(T_{mn}(k)) =  
  \sum_{B\subseteq \cE_{m,n}} 
     \frac{1}{[\cE_{m,n}-B]}
     \sum_{B\subseteq A,\;\mu(A)<k} (-1)^{|A-B|}.
\end{equation}

Hence, we have reduced the problem of computing $\rE(T_{mn}(k))$ to a sum 
over edge sets $B$ where each term is a fraction. The denominators are simple,
but the numerators are sums over supersets $A\supseteq B$ that do not
contain a $k$-assignment:
\begin{equation}\label{eq:ass-inner}
 S_{m,n,k}(B) := \sum_{B\subseteq A,\;
\mu(A)<k} (-1)^{|A-B|},
\end{equation} 
which turns out to be much more complicated than in the previous problem.
Note that $S_{m,n,k}(\emptyset)=0$ if $k=0$ (no terms in the sum) and also
if $k>m$ or $k>n$ (toggle any edge to see that the terms cancel). 

An edge set $B$ can be identified with a subset of the entries 
of an $m\times n$ matrix. 
Note that the existence of a $k$-assignment is independent
of the order of rows and columns. We say that a proper subset $B$ 
is {\em tabloidal} if the rows
and columns can be reordered to make $B$ tableau shaped, i.e.\@ 
consisting of all entries above a lattice path. It turns out that
only tabloidal $B$ give nonzero $S_{m,n,k}(B)$.

\begin{lemma}\label{th:nontableau}
$S_{m,n,k}(B)\ne 0$ only for tabloidal edge sets $B$.
\end{lemma}
\begin{proof}
It is easy to see that $B$ is nontabloidal if and only if its
set of matrix entries contains the pattern
{\setlength{\unitlength}{0.2mm}
\begin{picture}(20,15)(0,5)
\put(0,0){\framebox(20,20)}
\put(0,15){\colorbox{light}{\makebox(0,0){}}}
\put(10,5){\colorbox{light}{\makebox(0,0){}}}
\end{picture}}
. For some of the supersets $A$ in (\ref{eq:ass-inner}), the
upper right square is filled, and these cancel the terms where
this square is removed (the toggling trick). We are left with
$A$ such that the inclusion of the upper right square would
complete a $k$-assignment (that clearly doesn't involve the
lower left square). Since $\mu(A)=k\!-\!1$, K\"onig's
theorem tells us that it can be covered by $ k\!-\!1$ rows and
columns. The upper right square isn't covered so the lower left
square must be covered. Therefore it can be toggled and so everything
cancels out.
\end{proof}

\begin{figure}[h]\label{fig:tableau}
\begin{picture}(120,100)(0,0)
\put(0,0){\framebox(80,70)}
\put(0,23){\colorbox{light}{\makebox(14,24){}}}
\put(0,53){\colorbox{light}{\makebox(34,4){}}}
\put(0,63){\colorbox{light}{\makebox(74,4){}}}
\dashline[-7]{2}(0,10)(60,70)
\dottedline{2}(20,0)(20,20)
\dottedline{2}(20,20)(40,20)
\dottedline{2}(40,20)(40,50)
\dottedline{2}(40,50)(80,50)
\put(37,-10){$n$}
\put(85,31){$m$}
\put(0,80){\line(1,0){60}}
\put(0,78){\line(0,1){4}}
\put(60,78){\line(0,1){4}}
\put(27,83){$k$}
\end{picture}
\begin{picture}(50,50)(0,0)
\put(0,0){\framebox(20,20)}
\dashline[-7]{2}(0,10)(10,20)
\put(6,-10){$n_1$}
\put(25,7){$m_1$}
\put(0,30){\line(1,0){10}}
\put(0,28){\line(0,1){4}}
\put(10,28){\line(0,1){4}}
\put(1,35){$k_1$}
\end{picture}
\begin{picture}(50,60)(0,0)
\put(0,0){\framebox(20,30)}
\dashline[-7]{2}(0,10)(20,30)
\put(6,-10){$n_2$}
\put(25,12){$m_2$}
\put(0,40){\line(1,0){20}}
\put(0,38){\line(0,1){4}}
\put(20,38){\line(0,1){4}}
\put(5,45){$k_2$}
\end{picture}
\begin{picture}(60,50)(0,0)
\put(0,0){\framebox(40,10)}
\dashline[-7]{2}(0,0)(10,10)
\put(15,-10){$n_3$}
\put(45,3){$m_3$}
\put(0,20){\line(1,0){10}}
\put(0,18){\line(0,1){4}}
\put(10,18){\line(0,1){4}}
\put(0,27){$k_3$}
\end{picture}
\caption{The decomposition into mid rectangles.}
\end{figure}
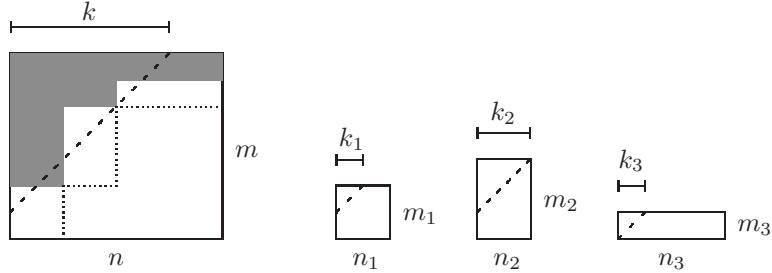

Not even all tabloidal edge sets $B$ contribute to the sum, only those in $\mbox{Tabloid}_k$,
that is such that the inner corners of the tableau
lie above the $k$-diagonal and the outer corners lie below
or on it (see Fig.\@ \ref{fig:tableau}). 
For such a tableau the numerator $S_{m,n,k}(B)$ is determined by the empty
mid rectangles
at the inner corners, defined by extending the vertical and horizontal lines.
Let the size of the $i$th rectangle be $m_i\times n_i$ and
note that $B\in\mbox{Tabloid}_k$ means that
the $k$-diagonal intersects the rectangle, defining a local $k_i$-diagonal.

\noindent
\begin{lemma}\label{th:decomposition}
With the notation of the figure, the following formula is true 
for any tabloidal $B\ne\emptyset$.
$$S_{m,n,k}(B)=(-1)^{I+1}\prod_{i=1}^IS_{m_i,n_i,k_i}(\emptyset).$$
Note that the product will be nonzero only if all $k_i>0$ (inner corners above
the $k$-diagonal) and all $k_i\le m_i,n_i$ (outer corners below or on
the $k$-diagonal). 
\end{lemma}
\begin{proof}
We must count extensions of $B$ not containing a $k$-assignment.
$S_{m,n,k}(B)$ is the number of even legal extensions minus the
number of odd legal extensions. 

First we show that it is enough to consider extensions obtained
by filling in entries in the mid rectangles only.
Let $U$ be the set of entries not in these rectangles nor in $B$.
To see that the contribution of all extensions with some fixed entry
in $U$ is zero we invoke Lemma~\ref{th:nontableau}. This argument
can be repeated for all entries in $U$ and,
by the inclusion-exclusion principle, the total contribution of extensions
with at least one entry in $U$ is
$$-\sum_{\emptyset\neq V\subseteq U}S_{m,n,k}(B+V)=0.$$

Then we show that an extension contains a $k$-assignment if and only if
the $i$th mid rectangle contains a $k_i$-assignment for all $i$.
Suppose that the latter statement is true. Then rows and columns can
be permuted such that the dashed $k$-diagonal is covered, thus producing
the $k$-assignment. Conversely, if the filled entries of the
$i$th mid rectangle can be covered by less than $k_i$ rows and columns,
then together with the rows above and the columns to the left we have
a covering of all filled entries by less than $k$ rows and columns.

Finally we note that the signed sum of extensions
not containing a $k$-assignment is equal to minus the corresponding
signed sum of extensions containing a $k$-assignment, since
$\sum (-1)^{|A|}=0$ where the sum is taken over all subsets $A$
of the union of the mid rectangles. Thus we get that
$$S_{m,n,k}(B)=(-1)^{I+1}\prod_{i=1}^IS_{m_i,n_i,k_i}(\emptyset).$$

\end{proof}

\noindent
Now, we are finally able to complete Theorem~\ref{th:time_assignment}.
\begin{propo}
In the formula
$$\rE(T_{mn}(k))=\sum_{B\in \mbox{\rm\footnotesize Tabloid}_k} 
{\frac{S_{m,n,k}(B)}{[\cE_{m,n}-B]}}$$
the numerators are given by
$$S_{m,n,k}(B)=(-1)^{I+1}\prod_{i=1}^I{\binom{m_i-1}{k_i-1}}
\binom{n_i-1}{k_i-1}$$
where $m_i$, $n_i$, $k_i$ refer to the $I$ mid rectangles at the inner
corners of the tableau $B$ (see Fig.\@ \ref{fig:tableau}).
\end{propo}
\begin{proof}

The only thing left to show is that
$$S_{m,n,k}(\emptyset)=\binom{m-1}{ k-1}\binom{n-1}{k-1}$$
and we will prove this by induction.
The statement is trivial for $m=n=1$ and for $k=1$.
Our induction assumption
will be that the statement is true for all values of $m,n,k$
smaller than the current ones (with at least one strict inequality).

We regroup the terms in $S_{m,n,k}(\emptyset)$
according to the number of filled entries in the first row.
Let $S_0$ denote the contribution from the legal extensions where the first row
is completely empty. By inclusion-exclusion we get
$$S_0=\sum_{V\subseteq \mbox{row}_1} S_{m,n,k}(V),$$
and by Lemma~\ref{th:decomposition}, all terms in the sum vanish except for
$V=\emptyset$ (the term we want to compute), $|V|=n$ (the whole row filled) and
$|V|=k-1$ (and there are $\binom{n}{k-1}$ such terms). It is obvious how
to reduce the cases of an empty row and of a filled row, so we get 
\begin{eqnarray*}
S_{m,n,k}(\emptyset)&=&S_{m-1,n,k}(\emptyset)
-\binom{n}{ k-1}(-1)S_{m-1,k-1,k-1}(\emptyset)S_{1,n-k+1,1}(\emptyset)
-S_{m-1,n,k-1}(\emptyset) \\
&=& \binom{m-2}{ k-1}\binom{n-1}{ k-1}
 + \binom{n}{ k-1}\binom{m-2}{ k-2}\binom{k-2}{ k-2}\binom{0}{ 0}\binom{n-k}{ 0} 
-\binom{m-2}{ k-2}\binom{n-1}{ k-2} \\ 
&=& \binom{m-1}{ k-1}\binom{n-1}{ k-1}.
\end{eqnarray*}
This completes the proof.
\end{proof}

\section{Conclusion}
\noindent
In this paper we succeed in finding formulas for $\rE(T_n(k))$, 
the expected time at which a spanning $k$-forest has appeared, for
$\rE(L_n(k))$, the expected length of this minimal $k$-forest, and for
$\rE(T_{mn}(k))$, the expected time at which a $k$-matching appears.
However, we fail to find a formula for
$\rE(L_{mn}(k))$, the expected length of a minimal $k$-matching.
Why is this?

The answer emerged in Sec.\ref{sec:forest}, in which we found that the
first three statistics shared an important property . 
Given the actual times $t_{ij}$
at which the edges $e_{ij}$ appeared, it is possible to express
$T_n(k)$, $L_n(k)$ and $T_{mn}(k)$ as sums of terms of the form
$\min(A)$, by which is meant the minimal time in the edge set $A$.

It is clear that many other statistics can be expressed analogously,
thereby giving rise to similar expectation formulas. However,
the minimal length assignment $L_{mn}(k)$ cannot be expressed by
subset minima in this way. We have been able to derive formulas
for $L_{mn}(k)$ for small $k$ (see \cite{EES}) and already $k=2$
gives a counterexample.

$$\rE(L_{mn}(2))=\frac{2}{[\cE_{m,n}]}+\sum_{i,j}\frac{[e_{ij}][\cE_{m,n}-e_{ij}]}
{[\cE_{m,n}][\cE_{m,n}-\mbox{row}_i][\cE_{m,n}-\mbox{col}_j]}.$$

It is easy to verify that this rational function cannot be rewritten
with linear denominators. Therefore, $L_{mn}(k)$ is not expressible
as a linear combination of minima.

\end{document}